\documentstyle[12pt]{article} 

\begin{document}

\begin{center}
\bf
\LARGE
Twenty Open Problems \\
in Enumeration of Matchings: \\
Progress Report \\
$ $ \\
\large
\rm
James Propp \\
September 2, 1998
\end{center}
\normalsize

\bigskip

\noindent
{\bf Problem 1:}
Two independent (and very different) solutions of this problem have been
found; one by Mihai Ciucu and Christian Krattenthaler,
and the other by Harald Helfgott.  Ciucu and Krattenthaler's preprint
``The number of centered lozenge tilings of a symmetric hexagon''
is available at
\begin{center}
{\tt http://radon.mat.univie.ac.at/People/} \\
{\tt kratt/artikel/fixrhomb.html}
\end{center}
and will appear in J.\ Combin.\ Theory Ser.\ A;
the authors compute more generally the number of rhombus tilings of a hexagon 
with sides $a,a,b,a,a,b$ that contain the central unit rhombus, where $a$
and $b$ must have opposite parity (the special case $a=2n-1$, $b=2n$ solves 
Problem 1).  The same generalization was obtained (in a different but 
equivalent form) by Harald Helfgott and Ira Gessel, using a completely 
different method; Helfgott and Gessel's preprint 
``Tilings of diamonds and hexagons with defects''
is available at
\begin{center}
{\tt http://www.cs.brandeis.edu/$\sim$ira/papers/enumtile.ps.gz}\ .
\end{center}
One might still try to look for a proof whose simplicity 
is comparable to that of the answer ``one-third''.

\medskip

\noindent
{\bf Problem 2:}
Christian Krattenthaler and Soichi Okada have evaluated the number of
rhombus tilings of an $(a,b+1,c,a+1,b,c+1)$-hexagon with the central
triangle removed; see
``The number of rhombus tilings of a `punctured' hexagon and the minor
summation formula'', Adv.\ Appl.\ Math.\ 21 (1998),
also available as
\begin{center}
{\tt http://radon.mat.univie.ac.at/People/} \\
{\tt kratt/artikel/punctured.html}\ .
\end{center}
By a different technique, Mihai Ciucu derived the
same formula in the case $a=b$.
Ira Gessel solved this problem independently using the determinant method;
his solution appears in the article with Helfgott
(see the last URL in the comments on Problem 1, above).

\medskip

\noindent
{\bf Problem 3:}
Theresia Eisenk\"olbl solved this problem.
What she does in fact is to compute the number of all rhombus tilings
of a hexagon with sides $a,b+3,c,a+3,b,c+3$ where an arbitrary triangle 
is removed from each of the ``long'' sides of the hexagon (not necessarily 
the triangle in the middle).
For the proof of her formula she uses nonintersecting lattice paths,
determinants, and the Desnanot-Jacobi determinant formula.
See ``Rhombus tilings of a hexagon with three missing border triangles'', 
available at the {\tt xxx} archives as article 
{\tt http://front.math.ucdavis.edu/math.CO/9712261}\ .

\medskip

\noindent
{\bf Problem 4:}
Theresia Eisenk\"olbl solved the first part of Problem 4 (no preprint
available at the time of this writing), and
Markus Fulmek and Christian Krattenthaler solved the second part.
Fulmek and Krattenthaler compute the number of rhombus tilings of a hexagon 
with sides $a,b,a,a,b,a$ (with $a$ and $b$ having the same parity)
that contain the rhombus touching the center of the hexagon and lying 
symmetric with respect to the symmetry axis that runs parallel to the 
sides of length $b$.
For the proof of their formula they compute Hankel determinants
featuring Bernoulli numbers, which they do by using facts about continued
fractions, orthogonal polynomials, and, in particular, continuous Hahn
polynomials. The special case $a=b$ solves the second part of Problem 4.
See Fulmek and Krattenthaler's articles ``The number of rhombus tilings of a 
symmetric hexagon which contain a fixed rhombus on the symmetry axis, I,'' 
Ann.\ Combin.\ 2 (1998), 19-40, also available as
\begin{center}
{\tt http://radon.mat.univie.ac.at/People/} \\
{\tt kratt/artikel/fixrhom2.html} \ ,
\end{center}
and ``The number of rhombus tilings of a symmetric hexagon which contain 
a fixed rhombus on the symmetry axis, II,'' 
available as
\begin{center}
{\tt http://radon.mat.univie.ac.at/People/} \\
{\tt kratt/artikel/fixrhom3.html}\ .
\end{center}

\medskip

\noindent
{\bf Problem 6:}
For Nicolau Saldanha's interpretation of the spectrum of $K K^*$, see
his preprint
``Generalized Kasteleyn matrices and their singular values,'' 
available as
\begin{center}
{\tt http://www.umpa.ens-lyon.fr/$\sim$nsaldanh/kk.ps.gz}\ .
\end{center}
Horst Sachs says that $K K^*$ may have some significance 
in the chemistry of polycyclic hydrocarbons (so-called
benzenoids) and related compounds as a useful approximate
measure of ``degree of aromaticity''.

\medskip

\noindent
{\bf Problem 8:}
Harald Helfgott has solved this problem.

\medskip

\noindent
{\bf Problem 9:}
This was already solved when I posed the problem; it is
a special case of Theorem 4.1 in Ciucu's paper ``Enumeration 
of perfect matchings in graphs with reflective symmetry'',
available as
\begin{center}
{\tt file://ftp.math.lsa.umich.edu/pub/ciucu/matchsymm.ps.Z}\ .
\end{center}

\medskip

\noindent
{\bf Problem 10:}
This problem was solved independently three times:
by Harald Helfgott and Ira Gessel,
by Christian Krattenthaler,
and by Eric Kuo.
Gessel and Helfgott
solve a more general problem than Problem 10
(for details, see the URL given above for Problem 1).
Krattenthaler's preprint 
``Schur function identities and the number of perfect matchings 
of holey Aztec rectangles'' is available as
\begin{center}
{\tt http://radon.mat.univie.ac.at/People/} \\
{\tt kratt/artikel/holeyazt.html}\ ;
\end{center}
it gives several results concerning the enumeration of perfect matchings
of Aztec rectangles where (a suitable number of) collinear vertices are
removed, of which Problem 10 is just a special case.  
There is some overlap between the results of Helfgott and Gessel 
and the results of Krattenthaler.

\medskip

\noindent
{\bf Problem 14:}
Constantin Chiscanu 
found a polynomial bound on 
the number of domino tilings of the Aztec window 
of inner order $x$ and outer order $x+w$;
for details, see
\begin{center}
{\tt http://www-math.mit.edu/$\sim$propp/chiscanu.ps.gz}\ .
\end{center}
Doug Zare used the transfer-matrix method to show that 
the number of tilings is not just bounded by a polynomial,
but given by a polynomial, for each fixed $w$; see
\begin{center}
{\tt http://www-math.mit.edu/$\sim$propp/zare}\ .
\end{center}

\medskip

\noindent
{\bf Problem 15:}
Ben Wieland solved this problem.

\medskip

\noindent
{\bf Problem 16:}
Ben Wieland solved this problem, too.
 
\medskip

\noindent
{\bf Problem 18:}
For Horst Sachs' response to this problem, see 
\begin{center}
{\tt http://math.mit.edu/$\sim$propp/kekule}\ .
\end{center}

\medskip

\noindent
{\bf Problem 19:}
Laszlo Lovasz gave a simple proof of my (oral) conjecture that
the number of perfect matchings of the $n$-cube 
has the same parity as $n$ itself.
Consider the orbit of a particular matching of the $n$-cube
under the group generated by the $n$ standard reflections of the $n$-cube.
If all the edges are parallel (which can happen in exactly $n$ ways),
the orbit has size 1;
otherwise the size of the orbit is 
of the form $2^k$ (with $k \geq 1$) --- an even number.
The claim follows,
and similar albeit more complex reasoning
should allow one to compute the enumerating sequence
modulo any power of 2.
Meanwhile,
L.H.\ Clark, J.C.\ George, and T.D.\ Porter,
in their article ``On the Number of 1-Factors in the $n$-Cube''
(Proceedings of the Twenty-eighth Southeastern
International Conference on Combinatorics, Graph Theory and Computing,
Congr.\ Numer.\ 127 (1997), 67--69)
show that if one lets $f(n)$ denote the number of 1-factors in the $n$-cube,
then $$f(n)^{2^{1-n}} \sim n/e$$
as $n \rightarrow \infty$.

\end{document}